\documentclass[10pt]{amsart}
\usepackage[square,compress,comma, numbers]{natbib}
\usepackage[colorlinks=true, citecolor=blue, linkcolor=blue]{hyperref}
\usepackage{amsfonts}
\allowdisplaybreaks[4]
\usepackage{amssymb}
\usepackage{amsmath}
\usepackage{color}
\newcommand{\abs}[1]{\left\lvert #1 \right\rvert}

\def\E#1{\mathbb{E}\left \{#1 \right\}}

\definecolor{c20}{rgb}{0.,0.7,0.}
\definecolor{c30}{rgb}{0.,0.,1.}
\definecolor{c40}{rgb}{1,0.1,0.7}
\definecolor{c50}{rgb}{1,0,0}
\definecolor{c60}{rgb}{1,0.9,0.1}
\definecolor{c70}{rgb}{0.50,1.00,0.00}

\def\N{\mathbb{N}}
\def\KD#1{{\textcolor{black}{#1}}}
\numberwithin{equation}{section}
\newtheorem{theo}{Theorem}[section]
\newtheorem{sat}[theo]{Proposition}
\newtheorem{de}[theo]{Definition}
\newtheorem{lem}{Lemma}[section]

\newtheorem{korr}[theo]{Corollary}
\newtheorem{remark}[theo]{Remark}
\newtheorem{remarks}[theo]{Remarks}
\numberwithin{equation}{section}

\newcommand{\prooftheo}[1]{ \textsc{Proof of Theorem} \ref{#1} }

\newcommand{\prooflem}[1]{\textsc{Proof of Lemma} \ref{#1}}

\newcommand{\pk}[1]{\mathbb{P} \left\{ #1 \right\} }

\newcommand{\QED}{\hfill $\Box$}
\newcommand{\COM}[1]{}
\def\IF{\infty}
\newcommand{\R}{\mathbb{R}}
\newcommand{\inr}{\in \R}

\newcommand{\BQN}{\begin{eqnarray}}
\newcommand{\EQN}{\end{eqnarray}}
\newcommand{\BQNY}{\begin{eqnarray*}}
\newcommand{\EQNY}{\end{eqnarray*}}

\def\polhk#1{\setbox0=\hbox{#1}{\ooalign{\hidewidth
\lower1.5ex\hbox{`}\hidewidth\crcr\unhbox0}}}
\def\pE#1{\textcolor{c20}{#1}}
\def\pE#1{#1}

\def\cL#1{\textcolor{c50}{#1}}
\def\cL#1{#1}

\def\tbb#1{\textcolor{blue}{#1}}
\def\tbb#1{#1}

\newcommand{\kb}[1]{\boldsymbol{#1}}
\newcommand{\vk}[1]{\kb{#1}}

\newcommand{\netheo}[1]{{Theorem \ref{#1}}}

\def\rw{\rightarrow}

\def\IF{\infty}

\def\Cov{\mathrm{Cov}}

\date{}

\def\oo{(1+o(1))}

\def\LT{\left}
\def\RT{\right}
\def\H{\mathcal{H}}

\def\rw{\rightarrow}

\def\TT{\mathcal{T}}

\def\Del{\triangle}

\def\vn{\varepsilon}
\def\Var{\text{Var}}

\def\FF{\widetilde{\mathcal{H}}}
\def\PP{\widetilde{\mathcal{P}}}
\def\H{\mathcal{H}}

\newcommand{\limit}[1]{\lim_{#1 \to \infty}}

\newcommand{\BS}{\begin{sat}}
\newcommand{\ES}{\end{sat}}
\newcommand{\BT}{\begin{theo}}
\newcommand{\ET}{\end{theo}}
\newcommand{\BK}{\begin{korr}}
\newcommand{\EK}{\end{korr}}

\newcommand{\BD}{\begin{de}}
\newcommand{\ED}{\end{de}}
\newcommand{\BIT}{\begin{itemize}}
\newcommand{\EIT}{\end{itemize}}
\newcommand{\BDI}{\begin{description}}
\newcommand{\EDI}{\end{description}}
\newcommand{\BRM}{\begin{remarks}}
\newcommand{\ERM}{\end{remarks}}
\newcommand{\BEL}{\begin{lem}}
\newcommand{\EEL}{\end{lem}}

\def\TT{\mathcal{T}}
\def\TT{\R }

\def\rw{\rightarrow}
\def\LT{\left}
\def\RT{\right}

\def\Var{\text{Var}}
\def\vn{\varepsilon}
\def\Cov{\mathrm{Cov}}

\def\oo{(1+o(1))}

\def\LT{\left}
\def\RT{\right}
\def\H{\mathcal{H}}

\def\rw{\rightarrow}

\def\TT{\mathcal{T}}

\def\Del{\triangle}

\def\vn{\varepsilon}
\def\Var{\text{Var}}

\def\oo{(1+o(1))}

\def\LT{\left}
\def\RT{\right}
\def\H{\mathcal{H}}

\def\rw{\rightarrow}

\def\TT{\mathcal{T}}

\def\Del{\triangle}

\def\vn{\varepsilon}
\def\Var{\text{Var}}

\def\FF{\widetilde{\mathcal{H}}}
\def\PP{\widetilde{\mathcal{P}}}
\def\H{\mathcal{H}}

\def\pE#1{\textcolor{c20}{#1}}
\def\pE#1{#1}

\def\cL#1{\textcolor{c50}{#1}}
\def\cL#1{#1}

\def\tbb#1{\textcolor{blue}{#1}}
\def\tbb#1{#1}
\begin{document}
\title{Parisian Ruin of the Brownian Motion Risk Model with Constant Force of Interest}
\author{Long Bai}
\address{Long Bai, Department of Actuarial Science, University of Lausanne, UNIL-Dorigny, 1015 Lausanne, Switzerland
}
\email{Long.Bai@unil.ch}
\author{Li Luo}
\address{Li Luo, Department of Actuarial Science, University of Lausanne, UNIL-Dorigny, 1015 Lausanne, Switzerland and
	School of Mathematical Sciences, Nankai University, Tianjin 300071, P.R. China}
\email{Li.Luo@unil.ch}
\bigskip
\maketitle
{\bf Abstract:} Let
$B(t), t\in \R$ be a standard Brownian motion. Define a risk process
\begin{align}\label{Rudef}
R_u^{\delta}(t)=e^{\delta t}\LT(u+c\int^{t}_{0}e^{-\delta s}d s-\sigma\int_{0}^{t}e^{-\delta s}d B(s)\RT),  t\geq0,
\end{align}
 where $u\geq 0$ is the initial reserve, $\delta\geq0$ is the force of interest, $c>0$ is the rate of premium and $\sigma>0$ is a volatility factor. In this contribution we obtain an approximation of the Parisian ruin probability
\begin{align*}
\mathcal{K}_S^{\delta}(u,T_u):=\pk{\inf_{t\in[0,S]} \sup_{s\in[t,t+T_u]} R_u^{\delta}(s)<0}, S\ge 0,
\end{align*}
as $u\rw\IF$ where $T_u$ is a bounded function. Further, we show that the Parisian ruin time of this risk process can be  approximated by an exponential random variable.
Our results are new even for the classical ruin probability and ruin time which correspond to $T_u\equiv0$ in the Parisian setting.\\

\def\TTT{\mathcal{T}}
\def\TT{\mathcal{T}}
\def\Z{\mathbb{Z}}
\def\HWD{\mathcal{H}_\alpha^\delta}
\def\phd{\mathcal{P}_{\alpha,\delta}^h}
\section{Introduction}
In a theoretical insurance model the surplus process $R_u(t)$ can be defined by
\BQNY
R_u(t)=u+ct-X(t), \ \ \ t\ge0,
\EQNY
see \cite{MR1458613}, where $u\geq 0$ is the initial reserve, $c>0$ is the rate of premium and $ X(t),t\geq0 $ denotes the aggregate claims process. More specifically, we assume that the aggregate claims process is a  Brownian motion, i.e., $X(t)=\sigma B(t),\ \sigma>0$.  Due to the nature of the financial market, we shall consider a more general surplus process including interest rate, see \cite{rolski2009stochastic}, called a risk reserve process with constant force of interest, i.e., $R_u^{\delta}(t),\ t\geq 0$, in \eqref{Rudef}.
See \cite{rolski2009stochastic, DHJ13a,HX2007} for more studies on risk models with force of interest.\\
During the time horizon $[0,S], S\in(0,\IF]$,  the classical ruin probability is defined as below
\BQN\label{fo3}
\psi_{S}^{\delta}(u):=\pk{\inf_{t\in[0,S]}R_u^{\delta}(t)<0},
\EQN
see \cite{MR1458613, husler1999extremes, MR2462285, dieker2005extremes}. In \cite{emanuel1975diffusion,harrison1977ruin} the exact formula of $\psi_{\IF}^{\delta}(u)$ for $\delta>0$ is shown to be
\BQNY
\psi_\infty^{\delta}(u)=
\frac{ \Psi\LT(\sqrt{\frac{2\delta}{\sigma^2}}u+\sqrt{\frac{2c^2}{\sigma^2\delta}} \RT)}{ \Psi\LT(\sqrt{\frac{2c^2}{\sigma^2\delta} }\RT)},\ \ u>0,
\EQNY
where $\Psi(x)=1- \Phi(x)$ with $\Phi(\cdot)$ the distribution function of an $\mathcal{N}(0,1)$ random variable.\\
For $\delta=0$, the exact value of $\psi_{\IF}^{0}(u)$ is well-known (cf. \cite{deelstra1994remarks}) with
\BQNY
\psi_\infty^{0}(u)=e^{-\frac{2cu}{\sigma^2}},\ \ u>0.
\EQNY
In the literature, there are no results for the classical ruin probability in the case of finite time horizon, i.e., $S\in(0,\IF)$.
For $S\in(0,\IF)$, with motivation from the recent contributions \cite{dkebicki2015parisian,debicki2015parisian} we shall investigate in this paper the Parisian ruin \KD{probability} over the time period $[0,S]$
defined as
\BQN\label{Parisian1}
\mathcal{K}_S^{\delta}(u,T_u):=\pk{\inf_{t\in[0,S]} \sup_{s\in[t,t+T_u]} R_u^{\delta}(s)<0},
\EQN
where $T_u\ge 0$ models the pre-specified time.  Our assumption on $T_u$ is that
$$\limit{u} T_uu^{2}= T\in[0,\IF) $$
and thus $\psi_S^{\delta}(u)$ is a special case of $\mathcal{K}_S^{\delta}(u,T_u)$ with $T_u\equiv0$.\\
Another quantity of interest is the conditional distribution of the ruin time for the surplus process $R_u^{\delta}(t)$. The classical ruin time, e.g., \cite{DHJ13a, HJ13, MR2462285}, is defined as
\BQN\label{fo21}
\tau(u)=\inf\{t>0: R_u^{\delta}(t)<0\}.
\EQN
Here as in \cite{dkebicki2015parisian} we define the Parisian ruin time of  the  risk process $R_u^{\delta}(t)$  by 
\BQN\label{eq:eta2}
\eta(u)=\inf\{t\ge T_u: t-\kappa_{t,u}\ge T_u, R_u^{\delta}(t)<0\},\ \ \ \text{with} \ \kappa_{t,u}=\sup\{s\in[0,t]: R_u^{\delta}(s)\ge0\},
\EQN
and $\tau(u)$ is a special case of $\eta(u)$ with $T_u\equiv 0$.\\
Brief organization of the rest of the paper: In Section 2 we first present our main results on the asymptotics of $\mathcal{K}_S^{\delta}(u,T_u)$ as $u\rw\IF$ and then
we display the approximation of the Parisian ruin time. All the  proofs are relegated to Section 3.

\section{Main results}
Before giving the main results, we shall introduce\COM{ the generalized Pickands constant as
\BQN\label{eq:GH}
\FF_{\alpha}(T)=\lim_{\lambda\to\IF}\frac{1}{\lambda}\FF_\alpha(\lambda,T),\ \ T\ge0,\  \alpha\in(0,2],
\EQN
with
\BQNY
\FF_\alpha(\lambda,T)=\E{\sup_{t\in[0,\lambda]}\inf_{ s\in[0,T]} e^{ \sqrt{2}B_{\alpha}(t-s)-\abs{t-s}^{\alpha}}}\in(0,\IF),\ \ \lambda,\ T\ge0,
\EQNY
where $B_\alpha(t)$ is a standard fractional Brownian motion (fBm) with Hurst index $\alpha/2\in(0,1]$.
Furthermore, we define} a generalized Piterbarg constant as
\BQN\label{eq:GP}
\PP(T)=\lim_{\lambda\to\IF} \PP(\lambda,T),\ \ T\ge0,
\EQN
where for $\lambda, T\ge0$
\BQNY
\PP(\lambda,T)=\E{\sup_{t\in[0,\lambda]}\inf_{ s\in[0,T]} e^{ \sqrt{2}B(t-s)-\abs{t-s}-(t-s)}}.
\EQNY
Note further that \COM{the  classical  Pickands constant $\H_\alpha$ equals $\FF_\alpha(0)$ and }the classical Piterbarg constant $\mathcal{P}_1^{1}[0,\IF)$ equals $\PP(0)$ and $\mathcal{P}_{1}^{1}[0,\IF)=2$, see \cite{MR1993262,GeneralPit16,MR3317852}. \\
\COM{It is known that$\mathcal{H}_{1}=1$ , $\mathcal{H}_{2}=1/\sqrt{\pi}$ ,$\mathcal{P}_{1}^{1}=1+\frac{1}{a}$ and $\mathcal{P}_{2}^{1}=\frac{1}{2}\left(1+\sqrt{1+\frac{1}{a}}\right)$, see \cite{MR1993262,GeneralPit16,MR3317852}.}
Through this paper $\sim$ means asymptotic equivalence when the argument tends to $0$ or $\IF$.
Recall that $\Psi(\cdot)$ denotes the tail distribution function of an $\mathcal{N}(0,1)$ random variable and  $\Psi(u)\sim\frac{1}{\sqrt{2\pi}u}e^{-\frac{u^2}{2}}, u\rw\IF$.
\COM{\BT\label{pan1}
For any $T\in(0,\infty)$ and $\delta>0$, we have
\BQN
\psi_{T}(u)\sim
\frac{\sigma}{\sqrt{\delta\pi}}\frac{\sqrt{1-e^{-2\delta T}}}{u+\frac{c}{\delta}(1-e^{-\delta T})}\exp\left(-\frac{\delta(u+\frac{c}{\delta}(1-e^{-\delta T}))^2}{\sigma^2(1-e^{-2\delta T})}\right),
\EQN
as $u\rw\IF$.
\ET}

\BT\label{paripan1}
For $ \delta>0,S>0$ and $\limit{u} T_uu^{2}= T\in[0,\IF) $, we have
\BQN
\mathcal{K}_S^{\delta}(u,T_u)\sim\PP(a T)\Psi\LT(\frac{\sqrt{2\delta}(u+\frac{c}{\delta}(1-e^{-\delta S}))}{\sigma\sqrt{1-e^{-2\delta S}}}\RT),\ u\rw\IF,
\EQN
where $a:=\frac{2\delta^2e^{-2\delta S}}{\sigma^2(1-e^{-2\delta S})^2}$.
\ET
\COM{The classical ruin time of the risk process $U(t)$ over $[0,T]$  is a special case of the Parisian ruin time as $T_u=0$ by
\BQN\label{fo2}
\tau(u)=\inf\{t\in [0,T] \mid U(t)<0\}.
\EQN}
\begin{remarks}\label{remark1}
a) When $T_{u}\equiv0$,  $\mathcal{K}_S^{\delta}(u,T_u)$ reduces to the classical ruin probability $\psi_S^{\delta}(u)$, and by Theorem \ref{paripan1} with $T=0$
\BQNY
\mathcal{K}_S^{\delta}(u,0)=\psi_S^{\delta}(u)\sim2\Psi\LT(\frac{\sqrt{2\delta}(u+\frac{c}{\delta}(1-e^{-\delta S}))}{\sigma\sqrt{1-e^{-2\delta S}}}\RT),\ u\rw\IF.
\EQNY
b) If $\delta=0$
\begin{align}\label{delta1}
\mathcal{K}_S^{0}(u,T_{u})&=\pk{\inf_{t\in[0,S]}\sup_{s\in[t,t+T_u]}\LT(u+cs-\sigma B(s)\RT)<0}\nonumber\\
&\sim\PP(b T)\Psi\LT(\frac{u+cS}{\sigma\sqrt{S}}\RT),\ u\rw\IF,
\end{align}
where $b:=\frac{1}{2\sigma^2S^2}$  and we used the result of Corollary 3.4 (ii) in \cite{debicki2015parisian}.\\
Further, if $\delta=0$ and $T_{u}\equiv0$, by (\ref{delta1}) with $T=0,$ we get the asymptotic result of the classical ruin probability
\BQN\label{delta2}
\psi_S^{0}(u)\sim2\Psi\LT(\frac{u+cS}{\sigma\sqrt{S}}\RT),\ u\rw\IF.
\EQN
In fact,  \cite{deelstra1994remarks} gave the exact result of $\psi_S^{0}(u),\ u>0$, i.e.,
\begin{align*}
\psi_S^{0}(u)&=\Psi\LT(\frac{u+cS}{\sigma\sqrt{S}}\RT)+e^{-\frac{2cu}{\sigma^2}}\Phi\LT(\frac{cS-u}{\sigma\sqrt{S}}\RT)\\
&\sim2\Psi\LT(\frac{u+cS}{\sigma\sqrt{S}}\RT), \ u\rw\IF,
\end{align*}
which follows from
\BQNY
\lim_{u\rightarrow\infty}\frac{e^{-\frac{2cu}{\sigma^2}}\Phi(\frac{cS-u}{\sigma\sqrt{S}})}{\Psi(\frac{u+cS}{\sigma\sqrt{S}})}
=\lim_{u\rightarrow\infty}\frac{-\frac{2c}{\sigma^2}e^{-\frac{2cu}{\sigma^2}}\Phi(\frac{cS-u}{\sigma\sqrt{S}})-\frac{1}{\sigma\sqrt{2\pi S}}
e^{-(\frac{u+cS}{\sigma\sqrt{S}})^2/2}}{-\frac{1}{\sigma\sqrt{2\pi S}}e^{-(\frac{u+cS}{\sigma\sqrt{S}})^2/2}}
=1.
\EQNY
\COM{In other words, the asymptotic result equals the exact value of $\psi_S(u)$ when u goes to infinity.}
\COM{c) Note that both $\FF_\alpha(\lambda,T)$ and $\PP^{b_1,b_2}_{\alpha}(\lambda,T)$ are well defined since
$$
\E{\exp\LT(\sup_{t\in[0,\lambda]}  \sqrt{2}B_{\alpha}(t)\RT)}<\IF,\ \ \forall \lambda\ge0,
$$
which follows directly from Piterbarg inequality (see Theorem 8.1 in \cite{Pit96} and its extension in Lemma 5.1 in \cite{KEP2015}).\\}
\end{remarks}
Our next result discusses the approximation of the conditional ruin time.
\BT\label{paripan2}
\COM{For any $T>0, x>0$, $\delta>0$ and $\tau(u)$ in \ref{fo21},
\BQN\label{fo4}
\pk{u^2(T-\tau(u))>x\mid \tau(u)\leq T}\sim\exp\left(-\frac{2\delta^2 e^{-2\delta T}}{\sigma^2(1-e^{-2\delta T})^2}x\right) ,  \ u\rightarrow\infty.
\EQN}
Let $\eta(u)$ satisfy (\ref{eq:eta2}), under the assumptions of Theorem \ref{paripan1}, we have for any $x>0$ and $\delta\geq0$,
\BQN\label{parifo4}
\pk{u^2(S+T_u-\eta(u))> x\big\lvert \eta(u)\leq S+T_u}\sim
\LT\{\begin{array}{ll}
\exp\LT(-a x\RT),&\ \text{if}\  \delta>0,\\
\exp\LT(-b x\RT), &\ \text{if}\  \delta=0,
\end{array}
\RT.\ u\rightarrow\IF,
\EQN
where $a:=\frac{2\delta^2e^{-2\delta S}}{\sigma^2(1-e^{-2\delta S})^2}$ and $b:=\frac{1}{2\sigma^2S^2}$. \COM{ That is, the random variable $\{u^2(S+T_u-\eta(u))\big|\eta(u)\leq S+T_u\}$ converges to a exponential random variable with parameter $\widehat{a}$ in distribution as $u\rw\IF$.}
\ET
\begin{remark}
If $T_{u}\equiv0$, then $\eta(u)=\tau(u)$ and by \netheo{paripan2}, we obtain as $u\rw\IF$
\BQNY
\pk{u^2(S-\tau(u))>x\mid \tau(u)\leq S}\sim
\LT\{\begin{array}{ll}
\exp\LT(-a x\RT),&\ \text{if}\  \delta>0,\\
\exp\LT(-b x\RT), &\ \text{if}\  \delta=0.
\end{array}
\RT.
\EQNY
\end{remark}
\section{Proofs}
Hereafter we assume that $ \mathbb{C}_i, i\in \N$ are some positive constants.
\COM{\prooftheo{pan1} For all $u$ large
\BQNY
\psi_T(u)&=&\pk{\sup_{t\in[0,T]}\LT(\sigma\int_{0}^{t}e^{-\delta s}d B(s)-c\int^{t}_{0}e^{-\delta s}d s\RT)>u}\\
&=&\pk{\sup_{t\in(0,T]}\overline{X}(t)\frac{f_{u}(T)}{f_{u}(t)}>f_{u}(T)},
\EQNY
with
\BQNY
X(t)=\sigma\int_{0}^{t}e^{-\delta s}d B(s),\
\overline{X}(t)=\frac{X(t)}{\sigma_{X}(t)}, \ \text{and}\
f_{u}(t)=\frac{u+c\int^{t}_{0}e^{-\delta s}d s}{\sigma_{X}(t)}=\frac{u+\frac{c}{\delta}(1-e^{-\delta t})}{\sigma_{X}(t)},\ t\in (0,T],
\EQNY
where $X(t)$ has variance function $V_{X}^{2}(t)=\E{X^{2}(t)}=\frac{\sigma^2}{2\delta}(1-e^{-2\delta t})>0, t\in(0,T].$\\
Obviously, $V_{X}^{2}(t)$ attains its maximum at $t=T$ with $V_{X}^{2}(T)=\frac{\sigma^2}{2\delta}(1-e^{-2\delta T})$.
We give the asymptotic expansion of the standard deviation function $V_{X}(t)$ at $T$, i.e.,
\BQN\label{fo51}
V_{X}(t)=\frac{\sigma}{\sqrt{2\delta}}\sqrt{1-e^{-2\delta T}}-\frac{\sigma}{\sqrt{2\delta}}\frac{e^{-\delta T}}{\sqrt{1-e^{-2\delta T}}}(e^{-\delta t}-e^{-\delta T})(1+o(1)), \ t\uparrow T.
\EQN
Moreover,
\BQNY
\Cov(X(t),X(s))=\frac{\sigma^2}{2\delta}(1-e^{-2\delta s}), \ t>s>0.
\EQNY
According to \eqref{fo51}, there exists a positive constant $a>0$ small enough such that
\BQN\label{1F2}
0\le e^{-\delta t}-e^{-\delta T}\leq C (V_{X}(T)-V_{X}(t))
\EQN
holds uniformly in $t\in[T-a,T]$ for some constant $C>0$.\\
For any $u>0$
\BQNY
\Pi(u)\leq\psi_T(u) \leq\Pi(u)+\pk{\sup_{t\in(0,T-a]}\overline{X}(t)\frac{f_{u}(T)}{f_{u}(t)}>f_{u}(T)},
\EQNY
where
\BQNY
\Pi(u)=\pk{\sup_{t\in[T-a,T]}\overline{X}(t)\frac{f_{u}(T)}{f_{u}(t)}>f_{u}(T)}.
\EQNY
Since
\BQN\label{fom11}
\frac{f_{u}(T)}{f_{u}(t)}=1-\frac{V_{X}(T)-V_{X}(t)}{V_{X}(T)}+\frac{ \frac{c}{\delta}(e^{-\delta t}-e^{-\delta T})}{u+\frac{c}{\delta}(1-e^{-\delta t})}\frac{V_{X}(t)}{V_{X}(T)},
\EQN
in view of \eqref{1F2}, for any $\varepsilon\in(0,1)$ and sufficiently large u
\BQNY
1-\frac{V_{X}(T)-V_{X}(t)}{V_{X}(T)}\leq\frac{f_{u}(T)}{f_{u}(t)}\leq 1-(1-\varepsilon)\frac{V_{X}(T)-V_{X}(t)}{V_{X}(T)}
\EQNY
holds uniformly in $t\in[T-a,T]$.  Consequently,
\BQNY
\pk{\sup_{t\in[T-a,T]}X_{0}(t)>f_{u}(T)}\leq \Pi(u)\leq \pk{\sup_{t\in[T-a,T]}X_{\varepsilon}(t)>f_{u}(T)},
\EQNY
where
\BQNY
X_{\varepsilon}(t):=\overline{X}(t)\left(1-(1-\varepsilon)\frac{V_{X}(T)-V_{X}(t)}{V_{X}(T)}\right),\ t\in[T-a,T],\ \varepsilon\in[0,1).
\EQNY
Assume that $\sigma_{X_{\varepsilon}}^{2}(t)$  is the variance of $X_{\varepsilon}(t)$, then
\BQNY
\sigma_{X_{\varepsilon}}^{2}(t)=\left(1-(1-\varepsilon)\frac{V_{X}(T)-V_{X}(t)}{V_{X}(T)}\right)^2
\EQNY
attains its unique maximum over $[T-a,T]$ at T with $\sigma_{X_{\varepsilon}}^{2}(T)=1$. Then by \eqref{fo51}
\BQN\label{Fom4}
\sigma_{X_{\varepsilon}}(t)=1-(1-\varepsilon)\frac{V_{X}(T)-V_{X}(t)}{V_{X}(T)}
=1-(1-\varepsilon)\frac{\delta e^{-2\delta T}}{1-e^{-2\delta T}}(T-t)(1+o(1)),  \ t\rightarrow T.
\EQN
In addition,
\BQNY
&&\Cov\left(\frac{X_{\varepsilon}(s)}{\sigma_{X_{\varepsilon}}(s)},\frac{X_{\varepsilon}(t)}{\sigma_{X_{\varepsilon}}(t)}\right)=\Cov(\overline{X}(t),\overline{X}(s))\\
&&\quad\quad=\sqrt{\frac{1-e^{-2\delta s}}{1-e^{-2\delta t}}}=1+\left(\sqrt{\frac{1-e^{-2\delta s}}{1-e^{-2\delta t}}}-1\right)\\
&&\quad\quad=1+\frac{e^{-2\delta t}-e^{-2\delta s}}{\sqrt{1-e^{-2\delta t}}(\sqrt{1-e^{-2\delta t}}+\sqrt{1-e^{-2\delta s}})}\\
&&\quad\quad=1-\frac{\delta e^{-2\delta T}|t-s|}{1-e^{-2\delta T}}(1+o(1)),\ \    s,t\rw T,
\EQNY
and furthermore, for some constant $C>0$
\BQNY
\E{(X_{\varepsilon}(t)-X_{\varepsilon}(s))^2}\leq 2C \abs{ t-s} , \quad s,t\in [T-a,T].
\EQNY
\COM{In order to obtain the constants in the form defined by the conditions of Theorem 10.1 in \cite{Pit20}, we need to scale the time, $t'=\frac{\delta e^{-2\delta T}t}{1-e^{-2\delta T}}$,  so that
\BQNY
\pk{\sup_{t\in[T-a,T]}X_{\varepsilon}(t)>f_{u}(T)}=\pk{\sup_{t'\in\LT[\frac{\delta e^{-2\delta T}}{1-e^{-2\delta T}}(T-a),\frac{\delta e^{-2\delta T}}{1-e^{-2\delta T}}T\RT]}X_{\varepsilon}\LT(\frac{1-e^{-2\delta T}}{\delta e^{-2\delta T}}t'\RT)>f_{u}(T)}.
\EQNY
Now the process on the right-hand side satisfies the conditions E1-E2 (cf. Chapter 10 in \cite{Pit20}) and in addition we note that $\alpha=\beta=1, a=(1-\varepsilon)$. Recalling $\Psi(u)\sim\frac{1}{\sqrt{2\pi}u}e^{-\frac{u^2}{2}}, \ u\rw \IF$,} Then by Theorem 8.2 of \cite{Pit96}
\BQNY
\pk{\sup_{t\in[T-a,T]}X_{\varepsilon}(t)>f_{u}(T)}\sim\Psi(f_{u}(T))\Psi(f_{u}(T))\mathcal{P}_{1}^{1-\varepsilon},\ u\rw\IF,
\EQNY
and letting $\varepsilon\rightarrow0$, we obtain the asymptotic upper bound of $\Pi(u)$.
Similarly, we drive the lower bound
\BQNY
\pk{\sup_{t\in[T-a,T]}X_{0}(t)>f_{u}(T)}\sim\Psi(f_{u}(T))\mathcal{P}_{1}^{1}, \ u\rw\IF.
\EQNY
Thus
\BQNY
\Pi(u)\sim\Psi(f_{u}(T))\mathcal{P}_{1}^{1}
\sim\frac{\sigma}{\sqrt{\delta\pi}}\frac{\sqrt{1-e^{-2\delta T}}}{u+\frac{c}{\delta}(1-e^{-\delta T})}\exp\left(-\frac{\delta(u+\frac{c}{\delta}(1-e^{-\delta T}))^2}{\sigma^2(1-e^{-2\delta T})}\right),\ u\rw\IF.
\EQNY

In order to complete the proof we need to show that
\BQN\label{F5}
\pk{ \sup_{t\in(0,T-a]}\overline{X}(t)\frac{f_{u}(T)}{f_{u}(t)}>f_{u}(T)}=o(\Pi(u)),  \  u\rightarrow\infty.
\EQN
In light of \eqref{fom11}, for all u sufficiently large
\BQNY
\sup_{t\in(0,T-a]}\Var\LT(\overline{X}(t)\frac{f_{u}(T)}{f_{u}(t)}\RT)\leq (\rho(a))^2<1,
\EQNY
where $\rho(a)$ is a positive function which exists due to  the continuity of $\LT\{\frac{f_{u}(T)}{f_{u}(t)},t\in(0,T]\RT\}$. Then by Borell inequality, we have
\BQNY
\pk{ \sup_{t\in(0,T-a]}\overline{X}(t)\frac{f_{u}(T)}{f_{u}(t)}>f_{u}(T)}\leq \exp{\left(-\frac{(f_{u}(T)-m)^2}{2\rho^2(a)}\right)}=o(\Pi(u)) ,\ \ u\rightarrow\infty,
\EQNY
where $m:=\E{\sup_{t\in (0,T-a]}\overline{X}(t)\frac{f_{u}(T)}{f_{u}(t)}}<\infty$.\\
Consequently, \eqref{F5} is established, and the proof is complete.
\QED}

\prooftheo{paripan1}
For $S>0$ and $u$ large enough
\begin{align*}
\mathcal{K}_S^{\delta}(u,T_u)&=\pk{\sup_{t\in[0,S]} \inf_{s\in[t,t+T_{u}]}\LT(\sigma\int_{0}^{s}e^{-\delta z}d B(z)-c\int^{s}_{0}e^{-\delta z}d z\RT)>u}\\
&=\pk{\sup_{t\in[0,S]}\inf_{s\in[t,t+T_{u}]}\overline{X}(s)\frac{f_{u}(S)}{f_{u}(s)}>f_{u}(S)}\\
&=\pk{\sup_{t\in[0,S]}\inf_{s\in[t,t+T_{u}]}X_{u}(s)>f_{u}(S)},
\end{align*}
with
\BQNY
X(s)=\sigma\int_{0}^{s}e^{-\delta z}d B(z),\
\overline{X}(s)=\frac{X(s)}{\sigma_{X}(s)},\
f_{u}(s)=\frac{u+\frac{c}{\delta}(1-e^{-\delta s})}{\sigma_{X}(s)} \ \text{and}
\ X_{u}(s)=\overline{X}(s)\frac{f_{u}(S)}{f_{u}(s)},
\EQNY
where $\sigma_{X}^{2}(s)$ is the variance of $X(s)$ with $\sigma_{X}^{2}(s)=\frac{\sigma^2}{2\delta}(1-e^{-2\delta s})$.\\
Set $\rho(u)=\LT(\frac{\ln u}{u}\RT)^2$  and for any $\lambda>0$, Bonferroni inequality yields
\BQN\label{bb}
\Pi_0(u):=\pk{\sup_{t\in[S-\lambda u^{-2},S]}\inf_{s\in[t,t+T_{u}]}X_{u}(s)>f_{u}(S)}\leq\mathcal{K}_S^{\delta}(u,T_u)\leq\Pi_0(u)+\Pi_1(u)+\Pi_2(u),
\EQN
where
\BQNY
\Pi_1(u)=\pk{\sup_{t\in[0,S-\rho(u)]}\inf_{s\in[t,t+T_{u}]}X_{u}(s)>f_{u}(S)},\ \
\Pi_2(u)=\pk{\sup_{t\in[S-\rho(u),S-\lambda u^{-2}]}\inf_{s\in[t,t+T_{u}]}X_{u}(s)>f_{u}(S)}.
\EQNY
First we give some upper bounds of $\Pi_1(u)$ and $\Pi_2(u)$ which finally show that
\BQN\label{up}
\Pi_1(u)+\Pi_2(u)=o\LT(\Pi_0(u)\RT), \ \ u\rw\IF.
\EQN
For all $u$ large
\begin{align*}
\E{(X_u(t_1)-X_u(t_2))^2}&=\E{\LT({X}(t_1)\frac{f_{u}(S)}{u+\frac{c}{\delta}\LT(1-e^{-\delta t_1}\RT)}-{X}(t_2)\frac{f_{u}(S)}{u+\frac{c}{\delta}\LT(1-e^{-\delta t_2}\RT)}\RT)^2}\\
&\leq\mathbb{C}_1\E{\LT(\int_{t_1}^{t_2}e^{-\delta z}d B(z)\RT)^2}+\mathbb{C}_2\LT(\frac{u+\frac{c}{\delta}(1-e^{-\delta S})}{u+\frac{c}{\delta}(1-e^{-\delta t_1})}-\frac{u+\frac{c}{\delta}(1-e^{-\delta S})}{u+\frac{c}{\delta}(1-e^{-\delta t_2})}\RT)^2\\
&\leq\mathbb{C}_3|t_1-t_2|,\ \ t_1<t_2,\   t_1,t_2\in(0,S].
\end{align*}
Moreover,
\BQNY
\sup_{t\in[0,S-\rho(u)]}\Var\LT(X_u(t)\RT)=\sup_{t\in[0,S-\rho(u)]}\LT(\frac{f_{u}(S)}{f_{u}(t)}\RT)^2
=\frac{f^2_u(S)}{f^2_u(S-\rho(u))},
\EQNY
where we use the fact that $f_u(t)$ is a decreasing function for $t\in[0,S]$ when $u$ large enough.
Therefore, by Theorem 8.1 in \cite{Pit96}, we obtain
\BQN\label{P1}
\Pi_1(u)\leq\pk{\sup_{t\in[0,S-\rho(u)]}X_{u}(t)>f_{u}(S)}\leq \mathbb{C}_4 u^2\Psi\LT(f_u(S-\rho(u))\RT),
\EQN
and direct calculation yields that
\begin{align*}
u^2\Psi\LT(f_u(S-\rho(u))\RT)&\leq \frac{u^2}{\sqrt{2\pi}f_u(S)}e^{-\frac{f_u^2(S)}{2}\LT(\frac{f_u^2(S-\rho(u))}{f_u^2(S)}-1\RT)}e^{-\frac{f_u^2(S)}{2}}\\
&\sim u^2 e^{-a(\ln u)^2}\Psi(f_u(S))=o\LT(\Psi(f_u(S))\RT), \ \ u\rw\IF,
\end{align*}
where $a=\frac{2\delta^2 e^{-2\delta S}}{\sigma^2(1-e^{-2\delta S})^2}$  and we use the fact that
\BQN\label{A1}
1-\frac{f_{u}(S)}{f_{u}(S-t)}\sim \frac{\delta e^{-2\delta S}}{1-e^{-2\delta S}}t,\ \ t\rw 0.
\EQN
Set
\BQNY
\Del_k=\LT[k\lambda u^{-2},(k+1)\lambda u^{-2}\RT],\  \cL{k\in \N}, \ \ \text{and}\ \
 N(u)=\LT\lfloor \lambda^{-1}\rho(u)u^2\RT\rfloor,
\EQNY
where $\lfloor\cdot\rfloor$ stands for the ceiling function, then
\begin{align}\label{eqp1}
\Pi_2(u)&\leq\pk{\sup_{t\in[S-\rho(u),S-\lambda u^{-2}]}X_{u}(t)>f_{u}(S)}\nonumber\\
         &=\pk{\sup_{t\in[\lambda u^{-2},\rho(u)]}X_{u}(S-t)>f_{u}(S)}\nonumber\\
         &\leq\sum_{k=1}^{N(u)}\pk{\sup_{t\in\Delta_k}X_{u}(S-t)>f_{u}(S)}\nonumber\\
         &\leq\sum_{k=1}^{N(u)}\pk{\sup_{t\in\Delta_0}\overline{X}(S-t)>f_{u}(S-k\lambda u^{-2})}\nonumber\\
         &\leq\sum_{k=1}^{N(u)}\pk{\sup_{t\in[0,\lambda]}\overline{X}(S-u^{-2}t)>f_{u}(S-k\lambda u^{-2})}.
\end{align}
Clearly,
\BQN\label{C1}
\inf_{1\leq k\leq N(u)}f_{u}(S-k\lambda u^{-2})\rw \IF, u\rw\IF,
\EQN
and for $t_1<t_2, \ t_1,t_2\in[0,S]$,
\BQNY
r_X(t_1,t_2):=\E{\overline{X}(t_1)\overline{X}(t_2)}=\sqrt{\frac{1-e^{-2\delta t_1}}{1-e^{-2\delta t_2}}}.
\EQNY
Further,
\BQN\label{C2}
&&\lim_{u\rw\IF}\sup_{1\leq k\leq N(u)}\underset{t_1,t_2\in[0,\lambda]}{\sup_{t_1\neq t_2,}}
\LT|f^2_{u}(S-k\lambda u^{-2})\frac{\Var{\LT(\overline{X}(S-u^{-2}t_1)-\overline{X}(S-u^{-2}t_2)\RT)}}{2a|t_1-t_2|}-1\RT|\nonumber\\
&&\quad\quad=\lim_{u\rw\IF}\sup_{1\leq k\leq N(u)}\underset{t_1,t_2\in[0,\lambda]}{\sup_{t_1\neq t_2,}}
\LT|f^2_{u}(S-k\lambda u^{-2})\frac{2-2r_X(S-u^{-2}t_1,S-u^{-2}t_2)}{2a|t_1-t_2|}-1\RT|\nonumber\\
&&\quad\quad=0,
\EQN
and
\BQN\label{C3}
&&\sup_{1\leq k\leq N(u)}\underset{t_1,t_2\in[0,\lambda]}{\sup_{|t_1-t_2|<\vn}}
f^2_{u}(S-k\lambda u^{-2})\E{\LT(\overline{X}(S-u^{-2}t_1)-\overline{X}(S-u^{-2}t_2)\RT)\overline{X}(S)}\nonumber\\
&&\quad\quad\leq \mathbb{C}_5 u^2\underset{t_1,t_2\in[0,\lambda]}{\sup_{|t_1-t_2|<\vn}}\LT|r_X(S-u^{-2}t_1,S)-r_X(S-u^{-2}t_2,S)\RT|\nonumber\\
&&\quad\quad\leq \mathbb{C}_6 u^2\underset{t_1,t_2\in[0,\lambda]}{\sup_{|t_1-t_2|<\vn}}\LT|\sqrt{1-e^{-2\delta(S-u^{-2}t_1)}}-\sqrt{1-e^{-2\delta(S-u^{-2}t_2)}}\RT|\nonumber\\
&&\quad\quad\leq \mathbb{C}_7 \underset{t_1,t_2\in[0,\lambda]}{\sup_{|t_1-t_2|<\vn}}|t_1-t_2|\rw 0 ,\ \ u\rw\IF, \ \vn\rw0.
\EQN
According to \eqref{C1}, \eqref{C2}, \eqref{C3} and Lemma 5.3 of \cite{KEP2015},  \eqref{eqp1} is followed by
\BQN\label{P2}
\Pi_2(u)\leq \mathbb{C}_8\lambda\sum_{k=1}^{N(u)}\Psi\LT(f_{u}(S-k\lambda u^{-2})\RT)
\leq\mathbb{C}_9 \Psi(f_u(S)) \lambda\sum_{k=1}^{\IF}e^{-\mathbb{C}_{10} k\lambda}= o\LT(\Psi(f_u(S))\RT),\  u\rw\IF, \ \lambda\rw\IF,
\EQN
where the last inequality follows from \eqref{A1}.\\
Next we give the asymptotic behavior of $\Pi_0(u)$ as $u\rw\IF$ based on an appropriate application of  the Appendix in \cite{debicki2015parisian}.
For any $\vn_1>0$ and  $u$ large enough
\begin{align*}
\Pi_0(u)&=\pk{\sup_{t\in[S-\lambda u^{-2},S]}\inf_{s\in[t,t+T_{u}]}X_{u}(s)>f_{u}(S)}\\
  &\leq\pk{\sup_{t\in[S-\lambda u^{-2},S]}\inf_{s\in[t,t+(1-\vn_1)Tu^{-2}]}X_{u}(s)>f_{u}(S)}\\
  &= \pk{\sup_{t\in[0,\lambda]}\inf_{s\in[0,(1-\vn_1)T]}X_{u}(S+u^{-2}s-u^{-2}t)>f_{u}(S)}\\
  &=\pk{\sup_{t\in[0,\lambda]}\inf_{s\in[0,(1-\vn_1)T]}Y_{u}(t,s)>f_{u}(S)}\\
  &=:\Pi^{+}_0(u)
\end{align*}
and
\BQNY
\Pi_0(u)\geq  \pk{\sup_{t\in[0,\lambda]}\inf_{s\in[0,(1+\vn_1)T]}Y_{u}(t,s)>f_{u}(S)}=:\Pi^{-}_0(u),
\EQNY
where $Y_u(t,s):=X_{u}(S+u^{-2}s-u^{-2}t), $ for $(t,s)\in[0,\lambda]\times[0,(1+\vn_1)T]$.\\
Since
\begin{align*}
\sigma_{Y_u}(t,s):&=\sqrt{\Var\LT(Y_u(t,s)\RT)}=\sqrt{\Var(X_u(S+u^{-2}s-u^{-2}t))}
=\frac{f_u(S)}{f_u(S+u^{-2}s-u^{-2}t)}
\end{align*}
and $\eqref{A1}$,
 there exists $d(t,s)=\frac{\delta e^{-2\delta S}}{1-e^{-2\delta S}}(t-s)$ such that
\BQN\label{variance1}
\lim_{u\rightarrow\infty}\sup_{(t,s)\in[0,\lambda]\times[0,(1+\vn_1)T]}\abs{u^2(1-\sigma_{Y_u}(t,s))-d(t,s)}=0.
\EQN
Moreover, for $(t_1,s_1),(t_2,s_2)\in[0,\lambda]\times[0,(1+\vn_1)T]$ and $s_1-t_1>s_2-t_2$,
\BQNY
&&\Var(Y_u(t_1,s_1)-Y_u(t_2,s_2))\\
&&\quad\quad=f_u^2(S)\E{\frac{X(S+u^{-2}s_1-u^{-2}t_1)}{u+\frac{c}{\delta}(1-e^{-\delta(S+u^{-2}s_1-u^{-2}t_1)})}-\frac{X(S+u^{-2}s_2-u^{-2}t_2)}{u+\frac{c}{\delta}(1-e^{-\delta(S+u^{-2}s_2-u^{-2}t_2)})}}^2\\
&&\quad\quad=f_u^2(S)(J_1(u)+J_2(u)+J_3(u)),
\EQNY
where
\begin{align*}
J_1(u)&=\E{\frac{X(S+u^{-2}s_1-u^{-2}t_1)-X(S+u^{-2}s_2-u^{-2}t_2)}{u+\frac{c}{\delta}(1-e^{-\delta(S+u^{-2}s_1-u^{-2}t_1)})}}^2,\\
J_2(u)&=2\frac{\frac{c}{\delta}(e^{-\delta(S+u^{-2}s_1-u^{-2}t_1)}-e^{-\delta(S+u^{-2}s_2-u^{-2}t_2)})}{(u+\frac{c}{\delta}(1-e^{-\delta(S+u^{-2}s_1-u^{-2}t_1)}))(u+\frac{c}{\delta}(1-e^{-\delta(S+u^{-2}s_2-u^{-2}t_2)}))}\\
&\quad\times\E{\LT(\frac{X(S+u^{-2}s_1-u^{-2}t_1)-X(S+u^{-2}s_2-u^{-2}t_2)}{u+\frac{c}{\delta}(1-e^{-\delta(S+u^{-2}s_1-u^{-2}t_1)})}\RT)X(S+u^{-2}s_2-u^{-2}t_2)}=0,\\
J_3(u)&=\LT(\frac{\frac{c}{\delta}(e^{-\delta(S+u^{-2}s_1-u^{-2}t_1)}-e^{-\delta(S+u^{-2}s_2-u^{-2}t_2)})}{(u+\frac{c}{\delta}(1-e^{-\delta(S+u^{-2}s_1-u^{-2}t_1)}))(u+\frac{c}{\delta}(1-e^{-\delta(S+u^{-2}s_2-u^{-2}t_2)}))}\RT)^2\E{X(S+u^{-2}s_2-u^{-2}t_2)}^2.
\end{align*}
Since
\begin{align*}\label{taylor}
\lim_{u\rightarrow\infty}u^2f_u^2(S)J_1(u)&=\lim_{u\rightarrow\infty}f_u^2(S)\E{X(S+u^{-2}s_1-u^{-2}t_1)-X(S+u^{-2}s_2-u^{-2}t_2)}^2\\
&=\lim_{u\rightarrow\infty}\frac{u^2}{\frac{\sigma^2}{2\delta}(1-e^{-2\delta S})}\frac{\sigma^2}{2\delta}(e^{-2\delta(S+u^{-2}s_2-u^{-2}t_2)}-e^{-2\delta(S+u^{-2}s_1-u^{-2}t_1)}) \\
&=\frac{2\delta e^{-2\delta S}}{1-e^{-2\delta S}}((s_1-s_2)-(t_1-t_2)) \\
&=\frac{2\delta e^{-2\delta S}}{1-e^{-2\delta S}}\Var\LT(B(s_1-t_1)-B(s_2-t_2)\RT),\\
\lim_{u\rightarrow\infty}u^2f_u^2(S)J_3(u)&\leq\lim_{u\rightarrow\infty}\mathbb{C}_{11}(e^{-\delta(S+u^{-2}s_1-u^{-2}t_1)}-e^{-\delta(S+u^{-2}s_2-u^{-2}t_2)})\E{X(S+u^{-2}s_2-u^{-2}t_2)}^2 =0,
\end{align*}
thus
\BQN\label{variance2}
\lim_{u\rightarrow\infty}u^2\Var(Y_u(t_1,s_1)-Y_u(t_2,s_2))=\frac{2\delta e^{-2\delta S}}{1-e^{-2\delta S}}\Var\LT(B(s_1-t_1)-B(s_2-t_2)\RT).
\EQN
Further, there exist some constant $G,u_0>0$, such that for any $u>u_0$
\BQN\label{hold}
u^2\Var(Y_u(t_1,s_1)-Y_u(t_2,s_2))\leq G(\abs{t_1-t_2}+\abs{s_1-s_2})
\EQN
holds uniformly with respect to $(t_1,s_1),(t_2,s_2)\in[0,\lambda]\times[0,(1+\vn_1)T]$. By \eqref{variance1}, \eqref{variance2}, \eqref{hold}, Lemma 5.1 in \cite{debicki2015parisian} and $\lim_{u\rightarrow\infty}f_u(S)/u=1/\sigma_X(S)$, we obtain
\BQN\label{parimain1}
\Pi^{-}_{0}(u)\sim\PP( a \lambda,  a (1+\vn_1)T)\Psi(f_u(S)), \ u\to\IF.
\EQN
Similarly
\BQNY
\Pi^{+}_{0}(u)\sim\PP( a \lambda,  a (1-\vn_1)T)\Psi(f_u(S)), \ u\to\IF.
\EQNY
Letting $\vn_1\rw 0$ and $\lambda\rw\IF$, we have
\BQNY
\Pi_0(u)\sim\PP( a T)\Psi(f_u(S)), u\rw\IF.
\EQNY
The above combined with \eqref{P1} and \eqref{P2} drives \eqref{up},
therefore by \eqref{bb} the proof is complete.\QED\\
\prooftheo{paripan2}
\COM{i) According to the definition of conditional probability, for any $x,u >0$, we have
\BQNY
\pk{u^2(T-\tau(u))>x \mid \tau(u)\leq T}=\frac{ \pk{ \sup_{t\in [0,T_{x}(u)]}W(t)> u}}{\pk{\sup_{t\in [0,T]}W(t)> u}},
\EQNY
where $T_{x}(u)=T-xu^{-2}$ and $W(t)=\sigma\int_{0}^{t}e^{-\delta s}d B(s)-c\int^{t}_{0}e^{-\delta s}ds$.\\
First, we shall investigate asymptotic of $ \pk{\sup_{t\in [0,T_{x}(u)]}W(t)> u}$, as $u\rightarrow\infty$.
For any $u$ large enough
\BQNY
\pk{\sup_{t\in[0,T_{x}(u)]}W(t)> u}=\pk{\sup_{t\in(0,T_{x}(u)]}\overline{X}(t)\frac{f_{u}(T_{x}(u))}{f_{u}(t)}>f_{u}(T_{x}(u))},
\EQNY
with
\BQNY
X(t)=\sigma\int_{0}^{t}e^{-\delta s}d B(s),\
\overline{X}(t)=\frac{X(t)}{V_{X}(t)}, \ \text{and}\
f_{u}(t)=\frac{u+c\int^{t}_{0}e^{-\delta s}d s}{V_{X}(t)}=\frac{u+\frac{c}{\delta}(1-e^{-\delta t})}{V_{X}(t)},
\EQNY
where $V_{X}^{2}(t)$ is the variance of $X(t)$. $V_{X}^{2}(t)=\frac{\sigma^2}{2\delta}(1-e^{-2\delta t})$ attains its maximum over $(0,T_{x}(u)]$ at $t=T_{x}(u)$ with $V_{X}^{2}(T_{x}(u))=\frac{\sigma^2}{2\delta}(1-e^{-2\delta T_{x}(u)})$.  Hence
\BQN\label{fo511}
\lim_{u\rightarrow\infty}\sup_{t\in(0,\delta(u)]}\abs{\frac{V_{X}(T_{x}(u))-V_{X}(T_{x}(u)-t)}{\frac{\sigma}{\sqrt{2\delta}}\frac{e^{-\delta T_{x}(u)}}{\sqrt{1-e^{-2\delta T_{x}(u)}}}(e^{-\delta (T_{x}(u)-t)}-e^{-\delta T_{x}(u)})}-1}=0,
\EQN
where $\delta(u):=(\frac{\ln u}{u})^2$.
\COM{Further,
\BQNY
\Cov(X(t),X(s))=\frac{\sigma^2}{2\delta}(1-e^{-2\delta s}), \ t>s>0.
\EQNY}
Therefore, for u sufficiently large, we have
\BQN\label{F21}
\sup_{t\in(0,\delta(u)]}\abs{\frac{e^{-\delta (T_{x}(u)-t)}-e^{-\delta T_{x}(u)}}{V_{X}(T_{x}(u)-t)-V_{X}(T_{x}(u))}}&\leq&\frac{\sqrt{2\delta}}{\sigma}\frac{\sqrt{1-e^{-2\delta T}}}{e^{-\delta T}}.
\EQN
Set
\BQNY
\Pi(u,x)&=&\pk{\sup_{t\in[T_{x}(u)-\delta(u),T_{x}(u)]}\overline{X}(t)\frac{f_{u}(T_{x}(u))}{f_{u}(t)}>f_{u}(T_{x}(u))}\\
&=&\pk{\sup_{t\in[0,\delta(u)]}X_u(t)>f_{u}(T_{x}(u))},
\EQNY
where $X_u(t)=\overline{X}(T_{x}(u)-t)\frac{f_{u}(T_{x}(u))}{f_{u}(T_{x}(u)-t)}$.  Then for any $u>0$
\BQNY
\Pi(u,x)\leq\pk{\sup_{t\in(0,T_{x}(u)]}W(t)>u} \leq\Pi(u,x)+\Pi_o(u,x),
\EQNY
where $\Pi_o(u,x):=\pk{\sup_{t\in[\delta(u),T_{x}(u))}X_u(t)>f_{u}(T_{x}(u))}$. Since
\BQN\label{fom111}
\frac{f_{u}(T_{x}(u))}{f_{u}(T_{x}(u)-t)}=1-\frac{V_{X}(T_{x}(u))-V_{X}(T_{x}(u)-t)}{V_{X}(T_{x}(u))}+\frac{ \frac{c}{\delta}(e^{-\delta (T_{x}(u)-t)}-e^{-\delta T_{x}(u)})}{u+\frac{c}{\delta}(1-e^{-\delta (T_{x}(u)-t)})}\frac{V_{X}(T_{x}(u)-t)}{V_{X}(T_{x}(u))},
\EQN
combining with \eqref{F21}, we have
\BQNY
\lim_{u\rightarrow\infty}\sup_{t\in(0,\delta(u)]}
\abs{\frac{1-\frac{f_{u}(T_{x}(u))}{f_{u}(T_{x}(u)-t)}}{\frac{V_{X}(T_{x}(u))-V_{X}(T_{x}(u)-t)}{V_{X}(T_{x}(u))}}-1}=0.
\EQNY
Denote $\sigma_{X_{u}}^2(t)$ the variance of $X_{u}$, and due to the uniformly continuous of function $\frac{\delta e^{-2\delta t}}{1-e^{-2\delta t}}$ over interval $[0,T]$, we have
\BQNY
\lim_{u\rightarrow\infty}\sup_{t\in(0,\delta(u)]}\abs{\sigma_{X_{u}}(t)-1+\frac{\delta e^{-2\delta T}}{1-e^{-2\delta T}}(T-t)}=0.
\EQNY
In addition,
\BQNY
\lim_{u\rightarrow\infty}\sup_{s,t\in(0,\delta(u)]}\abs{\Cov\LT(\frac{X_{u}(t)}{\sigma_{X_{u}}(t)},\frac{X_{u}(s)}{\sigma_{X_{u}}(s)}\RT)-1+\frac{\delta e^{-2\delta T}}{1-e^{-2\delta T}}\abs{t-s}}=0,
\EQNY
and for some constant $C>0$,
\BQNY
\E{(X_{u}(t)-X_{u}(s))^2}\leq 2C \abs{t-s} , \quad s,t\in (0,\delta(u)].
\EQNY
We can use Lemma 5.1 in the Appendix \cite{debicki2015parisian} to get that
\COM{
In order to obtain the constants in the form defined by the conditions of Theorem 10.1 in \cite{Pit20}, we need to scale the time, $t'=\frac{(1-\varepsilon_1)\delta e^{-2\delta T}t}{1-e^{-2\delta T}}$,  so that
\BQNY
\pk{\sup_{t\in[0,a]}X_{\varepsilon}(T_{x}(u)-t)>f_{u}(T_{x}(u))}=\pk{\sup_{t'\in\LT[0,\frac{(1-\varepsilon_1)\delta e^{-2\delta T}}{1-e^{-2\delta T}}a\RT]}X_{\varepsilon}\LT(\frac{1-e^{-2\delta T}}{(1-\varepsilon_1)\delta e^{-2\delta T}}(T_{x}(u)-t')\RT)>f_{u}(T_{x}(u))}.
\EQNY
Now the process on the right-hand side satisfies the conditions E1-E2 (cf. Chapter 10 in \cite{Pit20})and in addition we note that $\alpha=\beta=1, a=(1-\varepsilon)$. Recalling $\Psi(u)\sim\frac{1}{\sqrt{2\pi}u}e^{-\frac{u^2}{2}}, \ u\rw \IF$ and letting $\varepsilon_1\rightarrow0$ ,thus
\BQNY
\pk{\sup_{t\in[0,a]}X_{\varepsilon}(T_{x}(u)-t)>f_{u}(T_{x}(u))}&\sim&\Psi(f_{u}(T_{x}(u)))\lim_{T\rightarrow\infty}\E{\exp\left(\max_{t\in[0,T]}\sqrt{2}B_{1}(t)-(1+(1-\varepsilon))t\right)}\\
&=&\Psi(f_{u}(T_{x}(u)))\mathcal{P}_{1}^{1-\varepsilon},
\EQNY
and letting $ \varepsilon\rightarrow0$, we obtain the asymptotic upper bound for $\Pi(u,x)$ on $ [T-a,T_{x}(u)]$.
Similarly,
\BQNY
\pk{\sup_{t\in[T-a,T_{x}(u)]}X_{0,u}(t)>f_{u}(T_{x}(u))}\sim\Psi(f_{u}(T_{x}(u)))\mathcal{P}_{1}^{1}, \ u\rw\IF.
\EQNY
 So}
\BQNY
\Pi(u,x)\sim\Psi(f_{u}(T_{x}(u)))\mathcal{P}_{1}^{1}
\sim\frac{\sigma}{\sqrt{\delta\pi}}\frac{\sqrt{1-e^{-2\delta T_{x}(u)}}}{u+\frac{c}{\delta}(1-e^{-\delta T_{x}(u)})}\exp\left(-\frac{\delta(u+\frac{c}{\delta}(1-e^{-\delta T_{x}(u)}))^2}{\sigma^2(1-e^{-2\delta T_{x}(u)})}\right),\ u\rw\IF.
\EQNY
In order to complete the proof we need to show further that
\BQN\label{F511}
\Pi_o(u,x)=o(\Pi(u,x)),  \  u\rightarrow\infty.
\EQN
\COM{According to the asymptotic expansion of the standard deviation function $V_{X}(t)$, there exists small enough $\theta_{0}$, some constant $C>0$ and for all u sufficiently large, such that
\BQN\label{uniform1}
\sup_{t\in(T-\theta_{0},T_{x}(u)-(\frac{\ln u}{u})^2]}\abs{\frac{e^{-\delta t}-e^{-\delta T_{x}(u)}}{V_{X}(t)-V_{X}(T_{x}(u))}}&\leq&C.
\EQN}
Since $V_{X}(t)$ is increasing and continuous, then for all $u$ large and some some $\theta>0$
\BQNY
\sup_{t\in[0,T-\theta]}V_{X}(t)<V_{X}(T-\theta)<V_{X}(T_{x}(u))
\EQNY
and
\BQNY
\sup_{t\in[\theta_{0},T_x(u)]}\sigma_{X_u}^{2}(t)\leq (\rho(\theta_{0}))^2<1.
\EQNY
Clearly,
\BQNY
\Pi_o(u,x)&\le&\pk{\sup_{t\in[\theta_0,T_x(u))}  X_u(t) >f_{u}(T_{x}(u))}+\pk{\sup_{t\in[\delta(u),\theta_0]}X_u(t) >f_{u}(T_{x}(u))}\\
&=:&\Pi_1(u,x)+\Pi_2(u,x).
\EQNY
By Borell-TIS inequality (cf. \cite{})
\BQNY
\Pi_1(u,x)\leq\exp\LT(-\frac{\LT(f_{u}(T_{x}(u))-m\RT)^2}{2 \rho^2(\theta_0)}\RT)=o(\Pi(u,x)),
\EQNY
where $m:=\E{\sup_{t\in [\theta_0,T]}X_{u}(t)}<\infty$.\\
In view of \eqref{fom111}, for small enough $\varepsilon, \varepsilon_1>0$ and sufficiently large $u$
\BQNY
\sigma_{X_{u}}(t)\leq 1-(1-\varepsilon)\frac{\delta e^{-2\delta T}}{1-e^{-2\delta T}}\LT(\frac{\ln u}{u}\RT)^2, \ \ t\in[\delta(u),\theta_0].
\EQNY
In addition,
\BQNY
\Cov\left(\frac{X_{u}(s)}{\sigma_{X_{u}}(s)},\frac{X_{u}(t)}{\sigma_{X_{u}}(t)}\right)&=&\Cov(\overline{X}(t),\overline{X}(s))\\
&=&1+\frac{e^{-2\delta t}-e^{-2\delta s}}{\sqrt{1-e^{-2\delta t}}(\sqrt{1-e^{-2\delta t}}+\sqrt{1-e^{-2\delta s}})}\\
&\geq&1- (1+\varepsilon_1)\frac{\delta e^{-2\delta T}|t-s|}{1-e^{-2\delta T}},\ \    s,t\in [\delta(u),\theta_0],
\EQNY
and
\BQNY
&&\E{(X_{u}(t)-X_{u}(s))^2}\\
&&\quad\quad\leq 2\LT(1-\Cov\left(\frac{X_{u}(s)}{\sigma_{X_{u}}(s)},\frac{X_{u}(t)}{\sigma_{X_{u}}(t)}\right)\RT)+(\sigma_{X_{u}}(t)-\sigma_{X_{u}}(s))^2\\
&&\quad\quad\leq C \abs{t-s} , \quad s,t\in [0,\theta_{0}].
\EQNY
Therefore, we can apply Theorem 8.1 in \cite{Pit96}
\BQNY
\Pi_2(u,x)&\leq& C u\exp\LT(-\frac{f_{u}^{2}(T_{x}(u))}{2 (1-(1-\varepsilon)\frac{\delta e^{-2\delta T}}{1-e^{-2\delta T}}(\frac{\ln u}{u})^2)^2}\RT)\\
&\le&C u\exp\LT(-\frac{f_{u}^{2}(T_{x}(u))\LT(1+(1-\varepsilon)\frac{\delta e^{-2\delta T}}{1-e^{-2\delta T}}(\frac{\ln u }{u})^2\RT)^2}{2 }\RT)=o(\Pi(u,x)),
\EQNY
as $u\rw\IF$, where C is independent of $u$. Consequently, \eqref{F511} is established. And we have
\COM{By time change, for any $x,u>0$ ,
\BQNY
\pk{\sup_{t\in[0,T_{x}(u)]}W(t)>u}=\pk{\sup_{t\in [0,1]}W(T_{x}(u)t)>u},
\EQNY
we can easily calculate the variance function of $W(t)$ , which is given by $\frac{\sigma^2}{2\delta}(1-e^{-2\delta t})$. The variance function attains its maximum on the set [0,1] at unique point $t=1$.\\
Following the result of Theorem \ref{pan1}, we have}
\BQNY
\pk{\sup_{t\in [0,T_{x}(u)]}W(t)> u}=\frac{\sigma}{\sqrt{\delta\pi}}\frac{\sqrt{1-e^{-2\delta T_{x}(u)}}}{u+\frac{c}{\delta}(1-e^{-\delta T_{x}(u)})}\exp\left(-\frac{\delta(u+\frac{c}{\delta}(1-e^{-\delta T_{x}(u)}))^2}{\sigma^2(1-e^{-2\delta T_{x}(u)})}\right)(1+o(1)), \  \ u \rightarrow\infty .
\EQNY
Therefore, we conclude that for $x>0$
\BQNY
&&\pk{u^2(T-\tau(u))>x\mid \tau(u)\leq T} \\
&&\quad\quad\sim\exp\left(-\frac{\delta}{\sigma^2}\left(\left(\frac{u+\frac{c}{\delta}(1-e^{-\delta T_{x}(u)})}{\sqrt{1-e^{-2\delta T_{x}(u)}}}\right)^2-\left(\frac{u+\frac{c}{\delta}(1-e^{-\delta T})}{\sqrt{1-e^{-2\delta T}}}\right)^2\right)\right)\\
&&\quad\quad=\exp\left(-\frac{\delta h(u,T,T_{x}(u))}{\sigma^2(1-e^{-2\delta T_{x}(u)})(1-e^{-2\delta T})}\right),\ u\rw\IF,
\EQNY
where
\BQNY
h(u,T,T_{x}(u))&=&(1-e^{-2\delta T})\LT(u+\frac{c}{\delta}(1-e^{-\delta T_{x}(u)})\RT)^2-(1-e^{-2\delta T_{x}(u)})\LT(u+\frac{c}{\delta}(1-e^{-\delta T})\RT)^2\\
&=&(\frac{2uc}{\delta}+\frac{2c^2}{\delta^2})(e^{-\delta T}-e^{-\delta T_{x}(u)})+(\frac{2c^2}{\delta^2}+u^2+\frac{2uc}{\delta})(e^{-2\delta T_{x}(u)}-e^{-2\delta T})\\
&&+(\frac{2uc}{\delta}+\frac{2c^2}{\delta^2})(e^{-2\delta T-\delta T_{x}(u)}-e^{-2\delta T_{x}(u)-\delta T}).
\EQNY
Since
\BQNY
e^{-\delta T_{x}(u)}-e^{-\delta T}\sim\delta xe^{-\delta T}u^{-2} ,  \ \
e^{-2\delta T_{x}(u)}-e^{-2\delta T}\sim 2\delta xe^{-2\delta T}u^{-2},\  \ u\rightarrow\infty,
\EQNY
it follows that
\BQNY
\pk{u^2(T-\tau(u))>x\mid \tau(u)\leq T}\sim\exp\left(-\frac{2\delta ^2 e^{-2\delta T}}{\sigma^2(1-e^{-2\delta T})^2}x\right) ,\  \ u\rightarrow\infty.
\EQNY \\}
{\bf Case 1} $\delta>0$: According to the definition of conditional probability, for any $x,u >0$
\BQN\label{bbb}
&&\pk{u^2(S+T_u-\eta(u))>x \mid \eta(u)\leq S+T_u}\nonumber\\
&&\quad\quad=\frac{ \pk{ \sup_{t\in [0,S-xu^{-2}]}\inf_{s\in [t,t+T_u]}\LT(\sigma\int_{0}^{s}e^{-\delta z}d B(z)-c\int^{s}_{0}e^{-\delta z}dz\RT)> u}}{\pk{\sup_{t\in [0,S]}\inf_{s\in [t,t+T_u]}\LT(\sigma\int_{0}^{s}e^{-\delta z}d B(z)-c\int^{s}_{0}e^{-\delta z}dz\RT)> u}}.
\EQN
Using the same notation of $X(s),\ \overline{X}(s),\ f_u(s),\ X_{u}(s),\ \sigma_X(s)$ as in the proof of Theorem 2.1, we have for $u$ large enough
\BQNY
&&\pk{ \sup_{t\in [0,S-xu^{-2}]}\inf_{s\in [t,t+T_u]}\LT(\sigma\int_{0}^{s}e^{-\delta z}d B(z)-c\int^{s}_{0}e^{-\delta z}dz\RT)> u}\\
&&\quad\quad=\pk{\sup_{t\in[0,S-xu^{-2}]}\inf_{s\in[t,t+T_{u}]}\overline{X}(s)\frac{f_{u}(S)}{f_{u}(s)}>f_{u}(S)}\\
&&\quad\quad=\pk{\sup_{t\in[0,S-xu^{-2}]}\inf_{s\in[t,t+T_{u}]}X_{u}(s)>f_{u}(S)},
\EQNY
Set $\rho(u)=\LT(\frac{\ln u}{u}\RT)^2$. For any $\lambda>0$, Bonferroni inequality yields
\BQN\label{bb1}
\Pi_0^{*}(u)\leq\pk{\sup_{t\in[0,S-xu^{-2}]}\inf_{s\in[t,t+T_{u}]}X_{u}(s)>f_{u}(S)}\leq\Pi_0^{*}(u)+\Pi_1^{*}(u)+\Pi_2^{*}(u),
\EQN
where
\BQNY
&&\Pi_0^{*}(u)=\pk{\sup_{t\in[S-xu^{-2}-\lambda u^{-2},S-xu^{-2}]}\inf_{s\in[t,t+T_{u}]}X_{u}(s)>f_{u}(S)},\\
&&\Pi_1^{*}(u)=\pk{\sup_{t\in[0,S-\rho(u)]}\inf_{s\in[t,t+T_{u}]}X_{u}(s)>f_{u}(S)},\\
&&\Pi_2^{*}(u)=\pk{\sup_{t\in[S-\rho(u),S-xu^{-2}-\lambda u^{-2}]}\inf_{s\in[t,t+T_{u}]}X_{u}(s)>f_{u}(S)}.
\EQNY
By \eqref{P1} and \eqref{P2} in the proof of \netheo{paripan1}, we know
\BQN\label{bb2}
\Pi_1^{*}(u)=o\LT(\Psi(f_u(S))\RT) ,\ \ u\rw\IF,
\EQN
 and
 \BQN\label{bb3}
 \Pi_2^{*}(u)\leq \pk{\sup_{t\in[S-\rho(u),S-\lambda u^{-2}]}\inf_{s\in[t,t+T_{u}]}X_{u}(s)>f_{u}(S)}= o\LT(\Psi(f_u(S))\RT), \ u\rw\IF,\ \lambda\rw\IF.
 \EQN
 Next we give the asymptotic behavior of $\Pi_0^{*}(u)$ as $u\rw\IF$. For any $\vn_1>0$ and $u$ large enough
\begin{align*}
\Pi_0^{*}(u)&=\pk{\sup_{t\in[S-xu^{-2}-\lambda u^{-2},S-xu^{-2}]}\inf_{s\in[t,t+T_{u}]}\overline{X}(s)\frac{f_{u}(S)}{f_{u}(s)}>f_{u}(S)}\\
&=\pk{\sup_{t\in[S-xu^{-2}-\lambda u^{-2},S-xu^{-2}]}\inf_{s\in[t,t+T_{u}]}\overline{X}(s)\frac{f_{u}(S-xu^{-2})}{f_{u}(s)}>f_{u}(S-xu^{-2})}\\
  &\leq \pk{\sup_{t\in[S-xu^{-2}-\lambda u^{-2},S-xu^{-2}]}\inf_{s\in[t,t+(1-\vn_1)Tu^{-2}]}\overline{X}(s)\frac{f_{u}(S-xu^{-2})}{f_{u}(s)}>f_{u}(S-xu^{-2})}\\
  &= \pk{\sup_{t\in[0,\lambda]}\inf_{s\in[0,(1-\vn_1)T]}\overline{X}(S+u^{-2}s-u^{-2}t-u^{-2}x)\frac{f_{u}(S-xu^{-2})}{f_{u}(S+u^{-2}s-u^{-2}t-u^{-2}x)}>f_{u}(S-xu^{-2})}\\
  &= \pk{\sup_{t\in[0,\lambda]}\inf_{s\in[0,(1-\vn_1)T]}Y_{u}^{*}(t,s)>f_{u}(S-xu^{-2})}\\
  &=:\Pi^{*+}_0(u),
\end{align*}
and
\BQNY
\Pi_0^{*}(u)\geq  \pk{\sup_{t\in[0,\lambda]}\inf_{s\in[0,(1+\vn_1)T]}Y_{u}^{*}(t,s)>f_{u}(S-xu^{-2})}=:\Pi^{*-}_0(u),
\EQNY
where $Y_u^{*}(t,s):=\overline{X}(S+u^{-2}s-u^{-2}t-u^{-2}x)\frac{f_{u}(S-xu^{-2})}{f_{u}(S+u^{-2}s-u^{-2}t-u^{-2}x)}, \ (t,s)\in[0,\lambda]\times[0,(1+\vn_1)T]$ and
$\sigma^2_{Y_u^{*}}(t,s):=\Var(Y_u^{*}(t,s))=\LT(\frac{f_{u}(S-xu^{-2})}{f_{u}(S+u^{-2}s-u^{-2}t-u^{-2}x)}\RT)^2$.\\
Using the similar argumentation as \eqref{variance1} in the proof of \netheo{paripan1}, we have
\BQN\label{C4}
\lim_{u\rightarrow\infty}\sup_{(t,s)\in[0,\lambda]\times[0,(1+\vn_1)T]}\abs{u^2(1-\sigma_{Y_u^{*}}(t,s))-d(t,s)}=0,
\EQN
with $d(t,s)=\frac{\delta e^{-2\delta S}}{1-e^{-2\delta S}}(t-s)$. Moreover, \eqref{variance2},  \eqref{hold} still hold for $Y_u^{*}(t,s)$ and $ (t_1,s_1), (t_2,s_2)\in[0,\lambda]\times[0,(1+\vn_1)T]$.

By Lemma 5.1 in \cite{debicki2015parisian} and $\lim_{u\rightarrow\infty}f_u(S)/u=1/\sigma_X(S)$, we obtain
\BQNY
\Pi^{*-}_{0}(u)\sim\PP( a \lambda, a (1+\vn_1)T)\Psi(f_u(S-xu^{-2}))\sim e^{- a x}\PP( a \lambda,  a (1+\vn_1)T)\Psi(f_u(S)), \ u\to\IF.
\EQNY
Similarly,
\BQNY
\Pi^{*+}_{0}(u)\sim e^{-a x}\PP( a \lambda,  a (1-\vn_1)T)\Psi(f_u(S)), \ u\to\IF.
\EQNY
Letting $\vn_1\rw 0$ and $\lambda\rw\IF$, we have
\BQNY
\Pi_0^{*}(u)\sim e^{-a x}\PP( a T)\Psi(f_u(S)), u\rw\IF.
\EQNY
The above combined with \eqref{bb1}, \eqref{bb2} and \eqref{bb3} derives that
\BQNY
\pk{\sup_{t\in[0,S-xu^{-2}]}\inf_{s\in[t,t+T_{u}]}X_{u}(s)>f_{u}(S)}\sim e^{-a x}\PP( a T)\Psi(f_u(S)),\ u\rw\IF.
\EQNY
Thus, the claim follows by using the results of \netheo{paripan1} and \eqref{bbb}.\\
{\bf Case 2} $\delta=0$:
\BQNY
\pk{u^2(S+T_u-\eta(u))> x\big\lvert \eta(u)\leq S+T_u}=\frac{\pk{\sup_{t\in[0,S-xu^{-2}]}\inf_{s\in[t,t+T_u]}(\sigma B(s)-cs)>u}}{\pk{\sup_{t\in[0,S]}\inf_{s\in[t,t+T_u]}(\sigma B(s)-cs)>u}}.
\EQNY
For $u$ large enough
\BQNY
\pk{\sup_{t\in[0,S-xu^{-2}]}\inf_{s\in[t,t+T_u]}(\sigma B(s)-cs)>u}=\pk{\sup_{t\in[0,S-x u^{-2}]}\inf_{s\in[t,t+T_u]}\widetilde{X}_u(s)>f_u(S)},
\EQNY
with
\BQNY
X(s)=\sigma B(s),\ \overline{X}(s)=\frac{B(s)}{\sqrt{s}},\ f_u(s)=\frac{u+cs}{\sigma\sqrt{s}} \ \text{and}\ \ \widetilde{X}_u(s)=\overline{X}(s)\frac{f_u(S)}{f_u(s)}.
\EQNY
Set $\rho(u)=\LT(\frac{\ln u}{u}\RT)^2$.  For any $\lambda>0$, Bonferroni inequality yields
\BQN\label{bbb1}
\widetilde{\Pi}_0(u)\leq\pk{\sup_{t\in[0,S-xu^{-2}]}\inf_{s\in[t,t+T_{u}]}\widetilde{X}_{u}(s)>f_{u}(S)}\leq\widetilde{\Pi}_0(u)+\widetilde{\Pi}_1(u)+\widetilde{\Pi}_2(u),
\EQN
where
\BQNY
&&\widetilde{\Pi}_0(u)=\pk{\sup_{t\in[S-xu^{-2}-\lambda u^{-2},S-xu^{-2}]}\inf_{s\in[t,t+T_{u}]}\widetilde{X}_{u}(s)>f_{u}(S)},\\
&&\widetilde{\Pi}_1(u)=\pk{\sup_{t\in[0,S-\rho(u)]}\inf_{s\in[t,t+T_{u}]}\widetilde{X}_{u}(s)>f_{u}(S)},\\
&&\widetilde{\Pi}_2(u)=\pk{\sup_{t\in[S-\rho(u),S-xu^{-2}-\lambda u^{-2}]}\inf_{s\in[t,t+T_{u}]}\widetilde{X}_{u}(s)>f_{u}(S)}.
\EQNY
Notice that for $u$ large enough
\begin{align*}
\E{(\widetilde{X}_u(t_1)-\widetilde{X}_u(t_2))^2}&=\frac{1}{S}\E{\LT(\frac{u+cS}{u+ct_1}B(t_1)-\frac{u+cS}{u+ct_2}B(t_2)\RT)^2}\\
&\leq\mathbb{C}_{12}\E{\LT(B(t_1)-B(t_2)\RT)^2}+\mathbb{C}_{13}\LT(\frac{u+c S}{u+ c t_1}-\frac{u+ c S}{u+ c t_2} \RT)^2\\
&\leq\mathbb{C}_{14}|t_1-t_2|,\ \ t_1<t_2,\   t_1,t_2\in(0,S],
\end{align*}
and
\BQNY
\sup_{t\in[0,S-\rho(u)]}\Var\LT(\widetilde{X}_u(t)\RT)=\sup_{t\in[0,S-\rho(u)]}\LT(\frac{f_{u}(S)}{f_{u}(t)}\RT)^2
=\frac{f^2_u(S)}{f^2_u(S-\rho(u))},
\EQNY
where we use the fact that $f_u(t)$ is a decreasing function for $t\in[0,S]$ when $u$ large enough.

Moreover,
\BQNY
1-\frac{f_{u}(S)}{f_{u}(S-t)}\sim \frac{1}{2S}t,\ \ t\rw 0,
\EQNY
\BQNY
\inf_{1\leq k\leq N(u)}f_{u}(S-k\lambda u^{-2})\rw \IF, u\rw\IF,
\EQNY
and for $t_1<t_2, \ t_1,t_2\in[0,S]$,
\BQNY
r_{\widetilde{X}}(t_1,t_2):=\E{\overline{X}(t_1)\overline{X}(t_2)}=\sqrt{\frac{t_1}{t_2}}.
\EQNY
Then
\BQN\label{C21}
&&\lim_{u\rw\IF}\sup_{1\leq k\leq N(u)}\underset{t_1,t_2\in[0,\lambda]}{\sup_{t_1\neq t_2,}}
\LT|f^2_{u}(S-k\lambda u^{-2})\frac{\Var{\LT(\overline{X}(S-u^{-2}t_1)-\overline{X}(S-u^{-2}t_2)\RT)}}{2b|t_1-t_2|}-1\RT|\nonumber\\
&&\quad\quad=\lim_{u\rw\IF}\sup_{1\leq k\leq N(u)}\underset{t_1,t_2\in[0,\lambda]}{\sup_{t_1\neq t_2,}}
\LT|f^2_{u}(S-k\lambda u^{-2})\frac{2-2r_{\widetilde{X}}(S-u^{-2}t_1,S-u^{-2}t_2)}{2 b|t_1-t_2|}-1\RT|=0,
\EQN
where $b=\frac{1}{2\sigma^2S^2}$,
and
\BQN\label{C31}
&&\sup_{1\leq k\leq N(u)}\underset{t_1,t_2\in[0,\lambda]}{\sup_{|t_1-t_2|<\vn}}
f^2_{u}(S-k\lambda u^{-2})\E{\LT(\overline{X}(S-u^{-2}t_1)-\overline{X}(S-u^{-2}t_2)\RT)\overline{X}(S)}\nonumber\\
&&\quad\quad\leq \mathbb{C}_{15} u^2\underset{t_1,t_2\in[0,\lambda]}{\sup_{|t_1-t_2|<\vn}}\LT|r_{\widetilde{X}}(S-u^{-2}t_1,S)-r_{\widetilde{X}}(S-u^{-2}t_2,S)\RT|\nonumber\\
&&\quad\quad\leq \mathbb{C}_{16} u^2\underset{t_1,t_2\in[0,\lambda]}{\sup_{|t_1-t_2|<\vn}}\LT|\sqrt{S-u^{-2}t_1}-\sqrt{S-u^{-2}t_2}\RT|\nonumber\\
&&\quad\quad\leq \mathbb{C}_{17} \underset{t_1,t_2\in[0,\lambda]}{\sup_{|t_1-t_2|<\vn}}|t_1-t_2|\rw 0 ,\ \ u\rw\IF, \ \vn\rw0.
\EQN
By Theorem 8.1 in \cite{Pit96} and Lemma 5.3 in \cite{KEP2015},  using the similar argumentation as in the proof of \netheo{paripan1},
we derive
\BQN\label{oo}
\widetilde{\Pi}_1(u)+\widetilde{\Pi}_2(u)= o\LT(\Psi(f_u(S))\RT),\  u\rw\IF, \ \lambda\rw\IF.
\EQN
Next we give the asymptotic behavior of $\widetilde{\Pi}_0(u)$ as $u\rw\IF$. For any $\vn_1>0$ and  $u$ large enough
\begin{align*}
\widetilde{\Pi}_0(u)&=\pk{\sup_{t\in[S-xu^{-2}-\lambda u^{-2},S-xu^{-2}]}\inf_{s\in[t,t+T_{u}]}\overline{X}(s)\frac{f_u(S)}{f_u(s)}>f_{u}(S)}\\
&=\pk{\sup_{t\in[S-xu^{-2}-\lambda u^{-2},S-xu^{-2}]}\inf_{s\in[t,t+T_{u}]}\overline{X}(s)\frac{f_u(S-xu^{-2})}{f_u(s)}>f_{u}(S-xu^{-2})}\\
&\leq \pk{\sup_{t\in[0,\lambda]}\inf_{s\in[0,(1-\vn_1)T]}\widetilde{Y}_{u}(t,s)>f_{u}(S-xu^{-2})}\\
&=:\widetilde{\Pi}^{+}_0(u)
\end{align*}
and
\BQNY
\widetilde{\Pi}_0(u)\geq  \pk{\sup_{t\in[0,\lambda]}\inf_{s\in[0,(1+\vn_1)T]}\widetilde{Y}_{u}(t,s)>f_{u}(S-xu^{-2})}=:\widetilde{\Pi}^{-}_0(u),
\EQNY
where $\widetilde{Y}_u(t,s):=\overline{X}(S+u^{-2}s-u^{-2}t-u^{-2}x)\frac{f_u(S-xu^{-2})}{f_u(S+u^{-2}s-u^{-2}t-u^{-2}x)}, $ for $(t,s)\in[0,\lambda]\times[0,(1+\vn_1)T]$.\\
Using the similar argumentation as \eqref{variance1}, \eqref{variance2} and \eqref{hold} in the proof of \netheo{paripan1}, we obtain that
\BQN\label{C4}
\lim_{u\rightarrow\infty}\sup_{(t,s)\in[0,\lambda]\times[0,(1+\vn_1)T]}\abs{u^2(1-\sigma_{\widetilde{Y}_u}(t,s))-\widetilde{d}(t,s)}=0,
\EQN
with $\widetilde{d}(t,s)=\frac{1}{2 S}(t-s)$ and $\sigma_{\widetilde{Y}_u}(t,s):=\sqrt{\Var(\widetilde{Y}_u(t,s))}$,
\BQNY
\lim_{u\rightarrow\infty}u^2\Var(\widetilde{Y}_u(t_1,s_1)-\widetilde{Y}_u(t_2,s_2))=\frac{1}{S}\Var\LT(B(s_1-t_1)-B(s_2-t_2)\RT),
\EQNY
and for some constant $G$ and all $u$ large enough
\BQNY
u^2\Var(\widetilde{Y}_u(t_1,s_1)-\widetilde{Y}_u(t_2,s_2))\leq G(\abs{t_1-t_2}+\abs{s_1-s_2})
\EQNY
uniformly for $(t_1,s_1), (t_2,s_2)\in[0,\lambda]\times[0,(1+\vn_1)T]$.

By Lemma 5.1 in \cite{debicki2015parisian} and $\lim_{u\rightarrow\infty}f_u(S)/u=\frac{1}{\sigma\sqrt{S}}$, we obtain
\BQNY
\widetilde{\Pi}^{-}_{0}(u)\sim \PP( b \lambda,  b (1+\vn_1)T)\Psi(f_u(S-xu^{-2}))\sim e^{-b x}\PP( b \lambda,  b (1+\vn_1)T)\Psi(f_u(S)), \ u\to\IF.
\EQNY
Similarly,
\BQNY
\widetilde{\Pi}^{+}_{0}(u)\sim e^{-b x}\PP( b \lambda, b (1-\vn_1)T)\Psi(f_u(S)), \ u\to\IF.
\EQNY
Letting $\vn_1\rw 0$ and $\lambda\rw\IF$, we have
\BQNY
\widetilde{\Pi}_0(u)\sim e^{-b x}\PP(b T)\Psi(f_u(S)),\ u\rw\IF.
\EQNY
The above combined with \eqref{bbb1} and \eqref{oo} leads to
\BQNY
\pk{\sup_{t\in[0,S-xu^{-2}]}\inf_{s\in[t,t+T_{u}]}\widetilde{X}_{u}(s)>f_{u}(S_{x}(u))}\sim e^{-b x}\PP( b T)\Psi(f_u(S)),\ u\rw\IF.
\EQNY
Using the above asymptotic equality and b) of Remarks \ref{remark1}, we obtain the results.

\QED
\COM{
\section{Appendix}
The next theorem consists of the well known results about excursion probabilities for non-stationary Gaussian processes, see \cite{Pit20}.
\BT\label{pan5}
For $T\in(0,\IF)$, let $\{X(t),t\in[0,T]\}$ be a Gaussian process with continuous paths, variance function $\sigma^2(t)$ and correlation function $r(s,t)$. Further, we introduce the following assumptions:\\
E1 $\sigma(t)$ attains its maximum over $[0,T]$ at an unique point $t_0\in[0,T]$, and for some positive $a,\beta$,
\BQNY
\sigma(t)=1-a\mid t-t_{0}\mid^\beta(1+o(1)), \ t\rightarrow t_{0}.
\EQNY
E2 (Local stationarity) For some $\alpha\in(0,2]$,
\BQNY
r(s,t)=1-\mid t-s\mid^\alpha(1+o(1)),\ s\rightarrow t_{0},t\rightarrow t_{0}.
\EQNY
E3  (Regularity) For some positive $\gamma$,G and all $s,t\in[0,T]$
\BQNY
\E{X(t)-X(s)}^2\leq G\mid t-s\mid^\gamma.
\EQNY
Then we have\\
(i)  if $\beta>\alpha$
\BQNY
\pk{\max_{t\in[0,T]}X(t)>u}=\frac{2\mathcal{H}_{\alpha}\Gamma(1/\beta)}{a^{1/\beta}\beta}u^{2/\alpha-2/\beta}\Psi(u)(1+o(1))
\EQNY
as $u\rightarrow\infty$;\\
(ii) if $\beta=\alpha$,
\BQNY
\pk{\max_{t\in[0,T]}X(t)>u}=\mathcal{P}_{\alpha}^{a}\Psi(u)(1+o(1))
\EQNY
as $u\rightarrow\infty$, where
\BQNY
0<\mathcal{P}_{\alpha}^{a}&=&\lim_{T\rightarrow\infty}\mathcal{P}_{\alpha}^{a}(T)<\infty,\\
\mathcal{P}_{\alpha}^{a}(T)&=&\E{\exp{\left(\max_{t\in[-T,T]}\chi(t)-a\mid t\mid^\alpha\right)}},
\EQNY
(iii) if $\beta<\alpha$,
\BQNY
\pk{\max_{t\in[0,T]}X(t)>u}=\Psi(u)(1+o(1))
\EQNY
as $u\rightarrow\infty$.\\
In case $t_{0}=0$, (i) holds after dividing the right-hand side by 2; (ii) holds with
\BQNY
\mathcal{P}_{\alpha}^{0,a}(T)=\E{\exp{\left(\max_{t\in[0,T]}\chi(t)-a\mid t\mid ^\alpha\right)}}
\EQNY
instead of   $\mathcal{P}_{\alpha}^{a}(T)$; and (iii) does not change. The same is valid for $t_{0}=T$.
\ET
\def\sigxiu{\sigma_{\xi_u}}
\def\bD{ \mathbf{ D} }
Let $\mathbf{D}$ be a compact set in $\R^n$, $n\in \mathbb{N}$ and suppose without loss of generality that $\mathbf{0}\in \mathbf{D}$.
Further, let $\{\xi_u(\vk{t}), \vk{t}\in\mathbf{D}\}$, $u>0$
be a family of centered Gaussian random fields with a.s. continuous sample paths and variance function $\sigma_{\xi_u}^2(\vk{\cdot})$. Below $||\cdot||$ stands for the Euclidean norm in $\R^n$.  We assume that $\xi_u$ satisfies the following conditions:

{\bf C1:} 
 $\sigxiu(\vk{0})=1$ for all  $u$ large, and there exists some  \cL{bounded measurable} function $d(\vk{\cdot})$ on $\bD$  such that
$$\lim_{u\rw\IF}\sup_{\vk{t}\in\bD}
\abs{u^2(1-\sigxiu(\vk{t}))-d(\vk{t})}=0.$$

{\bf C2:} There exist some centered Gaussian random field $\{\eta(\vk{t}), \vk{t}\in\R^n\}$ with a.s. continuous sample paths, $\eta(\vk{0})=0$ and variance function $\sigma_\eta^2(\vk{\cdot})$ such that
$$
\lim_{u\rw\IF}
 u^2\Var(\xi_u(\vk{t})-\xi_u(\vk{s}))=2\Var(\eta(\vk{t})-\eta(\vk{s}))
 $$
  holds for all $\vk{s},\vk{t}\in  \bD$.

 {\bf C3:} There exist some  constants $G, \nu>0, u_0>0,$ such that, for any $u>u_0$  
\BQNY
u^2\Var(\xi_u(\vk{t})-\xi_u(\vk{s})) \le G \abs{ \abs{\vk{t}-\vk{s}}}^{\nu}
\EQNY
holds uniformly with respect to $\vk{t},\vk{s}\in  \bD$.

As in \cite{DebKo2013} let $F: C(\bD)\to \R$   be a continuous functional acting on $C(\bD)$, the space of continuous functions on the compact set $\bD$. Assume that:

{\bf F1:} $\abs{F(f)}\le \sup_{\vk{t}\in\bD} \pE{\abs{f(\vk{t})}}$ for any $f\in C(\bD)$.

{\bf F2:} $F(af+b)=aF(f)+b$ for any $f\in C(\bD)$ and $a>0,b\inr$.\\
For any \cL{bounded measurable} function $d(\vk{\cdot})$ on $\bD$ with $d(\vk{0})=0$ and  $F$ satisfying {\bf F1} we define a constant
\BQN\label{eq:GPitCons}
\H_{\eta,d}^F(\bD)=\E{\exp\LT(F\LT(\sqrt{2}\eta(\vk{t})-\sigma_\eta^2(\vk{t})-d(\vk{t})\RT)  \RT)}.
\EQN
Along the lines of the proof in \cite{DebKo2013} we get that $\H_{\eta,d}^F(\bD)\in(0,\IF)$.
\tbb{The following result generalizes} Lemma 6.1 in \cite{Pit96} and Lemma 1 in \cite{DebKo2013}.
\BEL\label{LemGP}
Let $\{\xi_u(\vk{t}), \vk{t}\in\bD\}$, $u>0$
be the family of centered Gaussian random fields defined as above satisfying {\bf C1-C3} with some function
 $d(\vk{\cdot})$ and some Gaussian random field $\eta$. Let $F: C(\bD)\to \R$   be a continuous functional such that  {\bf F1-F2} hold.
 Then, for any positive  measurable  function $g(\cdot)$   satisfying $\lim_{u\rw\IF}g(u)/u=a \in (0,\IF)$
\BQN\label{eq:lem1}
\pk{F(\xi_u)>g(u)}=\H_{a\eta,a^2d}^F(\bD)\Psi(g(u))\oo
\EQN
holds as $u\to\IF$, provided that $\pk{F(\xi_u)>g(u)}>0$ for all large $u$.
\EEL
\def\xiu{\xi_u}

\prooflem{LemGP} The proof is based on the classical approach rooted in the ideas of \cite{PicandsA, Pit96}. For all
$u>0$ large
\BQN\label{eq:lamS1}
\ \pk{F(\xi_u)>g(u)}=\Psi(g(u))\int_{\R} \exp\LT(w-\frac{w^2}{2(g(u))^2}\RT)\pk{F(\xi_u)>g(u)\Bigl| \xi_u(\vk{0})=g(u)-\frac{w}{g(u)}}dw.
\EQN
Let, for any $u>0, w\inr$,
$\zeta_u=\{\zeta_u(\vk{t})=g(u)(\xi_u(\vk{t})-g(u))+w,   \vk{t}\in\bD\}.$
Using {\bf F2} the conditional probability in the integrand of \eqref{eq:lamS1} can be written as
\BQNY
 \pk{F(\xi_u)>g(u)\Bigl| \xi_u(\vk{0})=g(u)-\frac{w}{g(u)}}=\pk{F(\chi_u)>w},
\EQNY
where $\chi_u=\zeta_u|\zeta_u(\vk{0})=0$.
Denote
\BQNY
R_{\xiu}(\vk{t},\vk{s})=\E{\xiu(\vk{t})\xiu(\vk{s})},\ \ \vk{s},\vk{t} \in\bD
\EQNY
to be the covariance function of $\xi_u$. We have that
the conditional random field
$\chi_u=
\LT\{\chi_u(\vk{t}) , \vk{t}\in\bD\RT\}
$
has the same finite-dimensional distributions as
\BQNY
\LT\{g(u)(\xiu(\vk{t})-R_{\xiu}(\vk{t},\vk{0})\xiu(\vk{0}))-(g(u))^2(1-R_{\xiu}(\vk{t},\vk{0}))+w(1-R_{\xiu}(\vk{t},\vk{0})), \vk{t}\in\bD\RT\}.
\EQNY

Therefore, the following convergence
\BQNY
\E{\chi_u(\vk{t})}=-(g(u))^2(1-R_{\xiu}(\vk{t},\vk{0}))+w(1-R_{\xiu}(\vk{t},\vk{0})) \to -a^2(\sigma_\eta^2(\vk{t})+d(\vk{t})), \  u\rw\IF 
\EQNY
holds, for any $w\in\R$, uniformly with respect to $\vk{t}\in\bD.$ Moreover, for any $\vk{t},\vk{s}\in\bD$ we have
\BQNY
\Var\Bigl( \chi_u(\vk{t})-\chi_u(\vk{s})\Bigr)&=&(g(u))^2\LT(\E{\Bigl( \xiu(\vk{t})-\xiu(\vk{s})\Bigr)^2}-\LT(R_{\xiu}(\vk{t},\vk{0})-R_{\xiu}(\vk{s},\vk{0})\RT)^2\RT)\\
& \to & 2a^2\Var(\eta(\vk{t})-\eta(\vk{s})), \  u\rw\IF.
\EQNY
Therefore, the finite-dimensional distributions of $\chi_u$ converge to those of $\widetilde{\eta}=\{\sqrt{2}a \eta(\vk{t})-\sigma_{a\eta}^2(\vk{t})-a^2d(\vk{t}), \vk{t}\in\bD\}$, whereas the tightness follows by Proposition 9.7 in \cite{Pit20}.
 The rest of the proof repeats  line-by-line
that of Lemma 1 in \cite{DebKo2013}.
\QED}

{\bf Acknowledgement}: Thanks to  Swiss National Science Foundation grant no.  200021-166274.
\bibliographystyle{plain}
\bibliography{EEECE}
\end{document}